\documentclass[11pt]{amsart}
\theoremstyle{remark}

\theoremstyle{definition}

\numberwithin{equation}{section}

\allowdisplaybreaks[4]

\DeclareMathOperator{\codim}{codim}
\DeclareMathOperator{\cd}{cd}

\DeclareMathOperator{\Pic}{Pic}

\DeclareMathOperator{\cl}{cl}
\DeclareMathOperator{\CH}{CH}
\DeclareMathOperator{\Hg}{Hg}

\DeclareMathOperator{\DIV}{div}

\def \BC{{\mathbb C}}
\def \BQ{{\mathbb Q}}

\def \BH{{\mathbb H}}

\def \BZ{{\mathbb Z}}

\def \P{{\mathbb P}}

\def \D{{\mathcal D}}

\def \D{{\mathcal D}}
\def \O{{\mathcal O}}

\def \isom{\cong}

\def \tensor{\mathop{\otimes}}
\def \wt#1{{\widetilde{#1}}}

\begin{document}

\title[Noether-Lefschetz for $K_1$]{Noether-Lefschetz for 
$K_1$ of a Surface, Revisited}
\author{Xi Chen}

\address{632 Central Academic Building\\ 
University of Alberta\\
Edmonton, Alberta T6G 2G1, CANADA} 
\email{xichen@math.ualberta.ca}

\author{James D. Lewis}

\address{632 Central Academic Building\\
University of Alberta\\
Edmonton, Alberta T6G 2G1, CANADA}
\email{lewisjd@gpu.srv.ualberta.ca}

\thanks{First author partially supported
by a startup grant from the University of Alberta. 
Second author partially supported
by a grant from the Natural Sciences and Engineering 
Research Council of Canada}

\date{December 13, 2002}
%\today

\abstract{Let $Z \subset \P^3$ be a general surface of
degree $d\geq 5$. Using a Lefschetz pencil
argument, we give a elementary new proof of the
vanishing of a regulator on $K_1(Z)$.}
\endabstract
\maketitle

\section{Statement of result}\label{SEC001} 
 
Let $Z$ be a smooth quasiprojective variety over $\BC$,
and for given nonnegative integers $k, m$,
let $\CH^k(Z,m)$ be the higher Chow group as introduced
in \cite{Blo1}. In \cite{Blo2}, Bloch constructs a cycle
class map into any suitable cohomology theory. In our setting,
the corresponding map is:
$$
\cl_{k,m} : \CH^k(Z,m) \to H^{2k-m}_{\D}(Z,\BQ(k)),
$$
where $H^{2k-m}_{\D}(Z,\BQ(k))$ is Deligne-Beilinson 
cohomology, which fits in a short exact sequence
$$
0 \to \frac{H^{2k-m-1}(Z,\BC)}{F^kH^{2k-m-1}(Z,\BC) + 
H^{2k-m-1}(Z,\BQ(k))}\to H^{2k-m}_{\D}(Z,\BQ(k))
$$
$$
\to F^kH^{2k-m}(Z,\BC) \bigcap H^{2k-m}(Z,\BQ(k)) \to 0.
$$
Our primary interest is when $Z$ is also complete, and
$m=1$. Thus one has the corresponding map:
$$
\cl_{k,1} : \CH^k(Z,1) \to \frac{H^{2k-2}(Z,\BC)}{F^kH^{2k-2}(Z,\BC) + 
H^{2k-2}(Z,\BQ(k))}.
$$
Let $\Hg^{k-1}(Z) := H^{2k-2}(Z,\BQ(k-1))\cap F^{k-1}H^{2k-2}(Z,\BC)$ 
be the Hodge group. Then one has an induced map
$$
\underline{\cl}_{k,1} : \CH^k(Z,1) \to 
\frac{H^{2k-2}(Z,\BC)}{F^kH^{2k-2}(Z,\BC) + 
\Hg^{k-1}(Z)\tensor \BC + H^{2k-2}(Z,\BQ(k))}.
$$
It is known that $\underline{\cl}_{k,1}$ is trivial 
for $Z$ a sufficiently general
complete intersection and of sufficiently high
multidegree. This is an consequence of the work
of Nori \cite{No}, together with a technique
similar to that given in \cite{G-S}. The argument
is presented in \cite{MS}. Further, it is noted
in \cite{MS}, based on an effective bound in
\cite{Pa}, that 
$$
\underline{\cl}_{2,1} : \CH^2(Z,1) \to 
\frac{H^{2}(Z,\BC)}{F^2H^{2}(Z,\BC) + 
\Hg^{1}(Z)\tensor \BC + H^{2}(Z,\BQ(2))},
$$
is trivial for sufficiently general surfaces
$Z\subset \P^3$ of degree $d\geq 5$. The method
of Nori involves passing to the universal family
of complete intersections of a given multidegree,
in a given projective space. A similar point
of view appears in \cite{Na}. In this paper, we give
an elementary and direct proof of the triviality of 
$\underline{\cl}_{2,1}$ for a general surface 
$Z\subset \P^3$ of degree $\geq 5$, by working 
with a Lefschetz pencil of degree $d\geq 5$ surfaces
in $\P^3$. Thus our main theorem is an elementary new
proof of the following: 
\bigskip

{\bf Main Theorem.} {\it For a sufficiently general surface
$Z\subset \P^3$ of degree $d\geq 5$, the map 
$\underline{\cl}_{2,1}$ is trivial.}
\bigskip

We remark that the theorem is trivially true,
without the generic hypothesis, if $\deg Z \leq 3$,
as $H^{2}(Z)$ is algebraic. From the works of Collino, Voisin,
S. M\"uller-Stach, {\it et al,}\ and more recently
the authors \cite{C-L}, it is false if $\deg Z = 4$.
Since our method requires only a Lefschetz pencil as
opposed to the universal family of surfaces of degree $d$
in $\P^3$, and that it provides a rather simple
proof of a counterexample of the Hodge-$\D$-conjecture
of Beilinson \cite{Bei1}, we believe that this approach 
has some merit. In particular, we believe that this argument
is potentially useful in other settings.

\section{Some definitions}\label{SEC002} 

(1) {\it Deligne cohomology.} We assume that the reader is familiar with Deligne
cohomology, such as can be found in \cite{Bei1} and
\cite{EV}. In the case of a smooth projective
variety $Z$, and if we put $\BQ(j) = \BQ(2\pi\sqrt{-1})^j$,
one introduces the Deligne complex
$$
\BQ(j)_{\D}: \quad
\BQ(j)\to {\O}_Z\to \Omega_Z^1\to\cdots\to \Omega_Z^{j-1},
$$
and defines $H_{\D}^i(Z,\BQ(j)) := \BH^i(\BQ(j)_{\D})$
(hypercohomology). This gives rise to a short exact sequence
$$
0 \to \frac{H^{i-1}(Z,\BC)}{F^jH^{i-1}(Z,\BC) + H^{i-1}(Z,\BQ(j))}
\to H_{\D}^i(Z,\BQ(j))
$$
$$ 
\to F^jH^i(Z,\BC)\bigcap H^i(Z,\BQ(j))
\to 0.
$$
A similar exact sequence holds quasiprojective $Z$ that
are not necessarily smooth. 
\bigskip

(2) {\it Higher Chow groups.} For a quasiprojective $Z$, 
the following abridged definition
of $\CH^k(Z,1)$ will suffice \cite{La} (cf. \cite{MS}). 
\bigskip

{\bf Definition.} $\CH^k(Z,1)$ is the homology of the
middle term in the complex
$$
\coprod_{\cd_ZY=k-2}K_2(\BC(Y)) \xrightarrow{\mathrm{Tame}}
\coprod_{\cd_ZY=k-1}K_1(\BC(Y)) \xrightarrow{\DIV} 
\coprod_{\cd_ZY=k}K_0(\BC(Y)),
$$
where we recall that $K_1(\bf F) = {\bf F}^\times$
and $K_0(\bf F) = \BZ$, for
a field ${\bf F}$, and Tame, div are
respectively the Tame symbol and divisor maps.
\bigskip

{\it Note:} For the most part, we will identify
$\CH^{k}(-,m)$ with $\CH^{k}(-,m)\tensor \BQ$, unless
there is a specific reason to work with 
$\CH^{k}(-,m)$ (and in which case the interpretation will be clear). 
\bigskip

(3) {\it Horizontal displacement.} Let  $h : W\to S$
be a proper smooth morphism of quasiprojective varieties over $\BC$,
where say for simplicity $\dim S = 1$,
with smooth projective fiber $W_t := h^{-1}(t)$.
Fix a reference point $t_0\in S$ and consider a disk $\Delta$
centered at $t_0$. It is well known that there
is a diffeomeomorphism $h^{-1}(\Delta) \approx \Delta\times W_{t_0}$.
Thus for a cohomology class $\gamma := \gamma_{t_0}\in H^\bullet(W_{t_0})$,
one can talk about its horizontal displacement 
$\gamma_t \in H^\bullet(W_t)$, for $t\in \Delta$ and more generally for
$t\in S$. Consider the Hodge decomposition
$H^\bullet(W_t,\BC) = \bigoplus_{p+q=\bullet}H^{p,q}(W_t)$,
$\gamma_t = \oplus_{p+q=\bullet}\gamma^{p,q}$.
We say that the Hodge $(p,q)$ components deform horizontally
if $\gamma_t^{p,q} = (\gamma^{p,q})_t$ for all $t\in \Delta$.
By analytic considerations of Hodge subbundles, this is
equivalent to saying that 
$\gamma_t^{p,q} = (\gamma^{p,q})_t$ for all $t\in S$. 

\section{Proof of the main theorem}\label{SEC003} 

Let $\{X_t\}_{t\in \P^1}$ be a Lefschetz pencil
of surfaces of degree
$d\ge 5$ in $\P^3$, i.e. the general fiber $X_t$ is
smooth, and each singular fiber has an ordinary double
point singularity. We will think of this pencil
in the form  $X\subset \P^3\times \P^1$, i.e.
where $X$ is the blowup of $\P^3$ along the base 
locus  $\cap_{t\in \P^1}X_t$.
Suppose that for a general $t\in \P^1$, 
the cycle class map $\cl_{2,1} : \CH^2(X_t, 1)\to
H^3_{\D}(X_t,{\BQ}(2))$ is nontrivial. 
We can assume that $X$ is defined over an algebraically
closed field $L$ of finite transcendence 
degree over $\BQ$, i.e. $X/_{\BC} = X_L\times \BC$. 
Let $\eta$ be the generic point of $\P^1_L$. For some
finite algebraic extension $K \supset L(\eta)$, and
via a suitable embedding $K \hookrightarrow \BC$, 
% namely for that general $t$, the  map 
% $K\hookrightarrow \BC$, induced by
% the evaluation map $L(\eta) \hookrightarrow \BC$ at $t$.
there is a class $\xi_K \in \CH^2(X_K := X_\eta\times K,1)$
such that $\cl_{2,1}(\xi_K) \ne 0$ in $H^3_{\D}(X_K(\BC),\BQ(2))$.
[The situation here is not unlike that found in \cite[p. 191]{Lew}.]
There is a smooth projective curve $\Gamma_L$ with function field 
$L(\Gamma) = K$. Then after a base change $Y = X\times_{\P^1} \Gamma$,
$\xi_K$ defines a cycle in $\xi\in \CH^2(Y_U, 1)$, where $U\subset \Gamma$
is a Zariski open subset of $\Gamma$ and $Y_U = \cup_{t\in U}
Y_t$. This uses the fact that
$$
\CH^2(X_K,1) = \CH^2(Y_{\tilde{\eta}},1) = 
\lim_{\buildrel \longrightarrow\over U}\CH^2(Y_U,1),
$$
where $Y_{\tilde{\eta}}$ is the generic fiber of $Y$ over $\Gamma_L$.
We want to spread $\xi$ to {\it all} of $\Gamma$. However, there is
obstruction preventing us to do it; rather we can extend it after
a suitable modification of $\xi$. That is, we will show that there
exists $\xi'\in \CH^2(Y,1)$ such that \(\cl_{2,1}(\xi_t) =
\cl_{2,1}(\xi_t')\) for every $t\in U$. 
Our main tool is the localization sequence\footnote{Strictly
speaking, we don't really need the localization sequence
in this paper. Rather, it is used out of convenience.}
\begin{equation}\label{E000}
\CH^2(Y, 1) \xrightarrow{} \CH^2(Y_U, 1) \xrightarrow{}
\CH^1(Y_B) \xrightarrow{} \CH^2(Y) \xrightarrow{} \CH^2(Y_U)
\xrightarrow{} 0
\end{equation}
over $\BQ$, where $B = \Gamma\backslash U$ and \(Y_B = \cup_{t\in B}
Y_t\).

Note that the map $\CH^1(Y_B) \to \CH^2(Y)$ might not be injective if \(|B|
> 1\), so there is obstruction to extend $\xi$ directly.

Let $H$ be a plane in $\P^3$ and $\pi^* H\subset Y$ be the pullback of
$H$ under the projection $\pi: Y\to \P^3$. Let $C_b = \pi^* H\cap Y_b$
for $b\in B$ and $C_B = \cup_{b\in B} C_b$. Let us first extend $\xi$
to $Y\backslash C_B$. We look at the localization sequence
\begin{equation}\label{E001}
\CH^2(Y\backslash C_B, 1) \xrightarrow{} \CH^2(Y_U, 1) \xrightarrow{}
\CH^1(Y_B\backslash C_B) \xrightarrow{} \CH^2(Y\backslash C_B)
\end{equation}
Note that
\begin{equation}\label{E002}
\CH^1(Y_B\backslash C_B) = {\mathop{\oplus}\limits_{b\in B}} 
\CH^1(Y_b\backslash C_b)
\end{equation}
We claim that $\CH^1(Y_t\backslash C_t) \tensor \BQ = 0$ for every
$t\in \Gamma$.

The classical Noether-Lefschetz theorem tells us that a general
surface of degree $d\ge 4$ in $\P^3$ has Picard rank $1$. This
statement was refined by Mark Green \cite{G} to the following. Let $M =
\P^N$ be the space parameterizing surfaces of degree $d$ in $\P^3$
and $M_2\subset M$ be the subset parameterizing surfaces with
Picard rank $\ge 2$. Then $\codim_M M_2 = d - 3$. So when \(d \ge
5\), $M_2$ has codimension at least $2$ in $M$ and a general pencil
will avoid this locus. Thus $\Pic(Y_t)\tensor \BQ = \BQ$ for every
$t\in \Gamma$. Note that $Y_t$ might be singular, i.e., $Y_t$ has
an ordinary double point. Since an ordinary double point is a
quotient singularity, every Weil divisor of $Y_t$ is
$\BQ$-Cartier. Therefore, \(\CH^1(Y_t)\tensor \BQ =
\Pic(Y_t)\tensor \BQ\). In any case, we have
\begin{equation}\label{E003}
\CH^1(Y_t)\tensor \BQ = \Pic(Y_t)\tensor \BQ = \Pic(\P^3)\tensor
\BQ = \BQ.
\end{equation}
Obviously, $\CH^1(Y_t)$ is generated by $C_t = \pi^* H\cap Y_t$ over
$\BQ$. Consequently,
\begin{equation}\label{E004}
\CH^1(Y_t\backslash C_t) \tensor \BQ = 0
\end{equation}
and there is no obstruction to extend $\xi$ to $Y\backslash C_B$. So
we may regard $\xi$ as a class in $\CH^2(Y\backslash C_B, 1)$ from now
on.

There might be obstruction to further extend $\xi$ to all of $Y$ by
the localization sequence
\begin{equation}\label{E005}
\CH^2(Y, 1) \xrightarrow{} \CH^2(Y\backslash C_B, 1) \xrightarrow{\phi}
\CH^0(C_B) \xrightarrow{\gamma} \CH^2(Y)
\end{equation}
where
\begin{equation}\label{E006}
\CH^0(C_B) = \mathop{\oplus}\limits_{b\in B} \CH^0(C_b) = \BQ^{\oplus \beta}
\end{equation}
with $\beta = |B|$.

Let $\xi = \sum_\alpha (f_\alpha, D_\alpha)$ where $D_\alpha$ is a divisor on
$Y\backslash C_B$ and $f_\alpha$ is a rational function on
$D_\alpha$. We have
\begin{equation}\label{E007}
\sum_\alpha \DIV(f_\alpha) = 0.
\end{equation}
Let $\overline{D}_\alpha$ be the closure of $D_\alpha$ in $Y$ and
$f_\alpha$ naturally extends to a rational function
$\overline{f}_\alpha$ on $\overline{D}_\alpha$. Let 
\(\overline{\xi} = \sum_\alpha (\overline{f}_\alpha,
\overline{D}_\alpha)\). We no longer have \eqref{E007}. Instead,
\begin{equation}\label{E008}
\sum_\alpha \DIV(\overline{f}_\alpha) = \sum_{b\in B} m_b C_b
\end{equation}
for some $m_b\in\BZ$. Actually, the RHS of \eqref{E008} is exactly the
image of $\xi$ under the map \(\phi: \CH^2(Y\backslash C_B, 1) \to
\CH^0(C_B)\) in \eqref{E005}, i.e.,
\begin{equation}\label{E009}
\phi(\xi) = \sum_{b\in B} m_b C_b.
\end{equation}
Note that $\phi(\xi)$ lies in the kernel of
$\gamma: \CH^0(C_B) \to \CH^2(Y)$ and there is a
natural map $\CH^0(C_B)\to \CH^1(\Gamma)$ via
\begin{equation}\label{E010}
\CH^0(C_B) \xrightarrow{\gamma} \CH^2(Y) \xrightarrow{}
\CH^3(\P^3\times\Gamma)
\xrightarrow{} \CH^1(\Gamma).
\end{equation}
Note that the map $\CH^{3}({\bf P}^{3}\times \Gamma)
\xrightarrow{} [\P^{1}]\tensor\CH^1(\Gamma) = \CH^1(\Gamma)$, comes from the 
projective bundle formula.
Of course, the map $\CH^0(C_B)\to \CH^1(\Gamma)$ simply sends $C_b$ to
$Nb$, where $N=d$. And $\phi(\xi)$ maps to zero under this map, 
i.e. the divisor $\sum m_b b$
is $N$-torsion in $\CH^1(\Gamma) = \Pic(\Gamma)$. 

Note that $\pi^* H$ is a fibration of curves over $\Gamma$.
% Flat pullback, viz. any dominant morphism of a
% variety to a nonsingular curve is flat.
So the fact $\sum m_b b$ is torsion in $\CH^1(\Gamma)$
implies that $\sum m_b C_b$ is $N$-torsion in \(\CH^1(\pi^*
H)\). Consequently, there exists a rational function $f_H$ on $\pi^* H$
such that
\begin{equation}\label{E011}
\DIV(f_H) = N\sum_{b\in B} m_b C_b.
\end{equation}
So we may simply modify $\overline{\xi}$ as follows
\begin{equation}\label{E012}
\xi' = \overline{\xi} - \frac{1}{N}(f_H, \pi^*H).
\end{equation}
Now $\xi'\in \CH^2(Y, 1)$ and $\underline{\cl}_{2,1}(\xi_t') = 
\underline{\cl}_{2,1}(\xi_t)$ for
all $t\in U$, where we recall that
$$
\underline{\cl}_{2,1} : \CH^2(Y_t,1) \to \frac{H^3_{\D}(Y_t,\BQ(2))}
{\Hg^1(Y_t)\tensor \big(\BC/\BQ(1)\big)}
$$
is the induced map. This is due to the fact that
the restrictions $f_H$ to $Y_t$ are obviously
constants. Thus we can now replace $\xi$ by $\xi'$. 
Next observe that even though $Y$ is complete,
it may be singular. It is worthwhile pointing
out that we can further pull back $\xi$ to
a desingularization $\wt{Y}$ of $Y$. More precisely,
\bigskip

{\it Claim.} There exists
$\wt{\xi}\in \CH^2(\wt{Y}, 1)$ such that $\wt{\xi}$ and $\xi$ agree on
the open set where $\wt{Y}$ and $Y$ are isomorphic. 
\bigskip

The usefulness of this claim is as follows. The (cohomological)
cycle class map $\cl_{2,1} : \CH^2(Y,1)\to H^3_{\D}(Y,\BQ(2))$
is only defined if $Y$ is smooth. Granting the existence
of this cycle class map, the remaining argument only requires
the completeness of $Y$.
There is a short exact sequence:
$$
0\to \frac{H^2(Y,\BC)}{F^2H^2(Y,\BC) + H^2(Y,\BQ(2))}\to
H^3_{\D}(Y,\BQ(2)) \to F^2\cap H^3(Y,\BQ(2))\to 0.
$$
But since $Y$ is complete, a weight argument gives
$F^2\cap H^3(Y,\BQ(2)) = 0$. Thus 
for $t\in U$, $\underline{\cl}_{2,1}(\xi_t)$ is given by
the restriction $\underline{\cl}_{2,1}(\xi)\big|_{Y_t}$, i.e. 
induced by the restriction
$$
\frac{H^2(Y,\BC)}{F^2H^2(Y,\BC) + H^2(Y,\BQ(2))} \to
\frac{H^2(Y_t,\BC)}{F^2H^2(Y_t,\BC) + H^2(Y_t,\BQ(2))}.
$$
Thus as $t\in U$ varies, the class $\underline{\cl}_{2,1}(\xi_t)$
varies by {\it horizontal} displacement; further, the
restriction $H^2(Y) \to H^2(Y_t)$ is a morphism of
mixed Hodge structures. Thus $\underline{\cl}_{2,1}(\xi_t)$
is induced by a class in $H^2(Y_t)$, whose Hodge $(p,q)$
components displace {\it horizontally}, i.e. preserving the given
Hodge type. But over the set where $\Gamma \to \P^1$ ramifies,
one can find open sets $\Delta_\Gamma \subset U\subset \Gamma$, $\Delta\subset
\P^1$, in the strong topology, such that $\Delta_\Gamma \simeq \Delta$.
Thus $\underline{\cl}_{2,1}(\xi_t) = 0$, by virtue of:
\bigskip

{\bf Lemma.} {\it Consider a Lefschetz pencil $\{Z_t\}_{t\in \P^1}$
of surfaces in $\P^3$ of degree $d\geq 1$, and let 
$U_0 \subset \P^1$ be the smooth set. Further, let 
$\Delta \subset U_0$ be a disk, and assume given
$\gamma_t \in H^2(Z_t,\BC)$, a horizontal displacement
of a class $\gamma$ for $t\in \Delta$. If the $(p,q)$
components of $\gamma_t$ also horizontally displace, then
$\gamma_t \in \Hg^1(Z_t)$.} 
\bigskip

{\it Proof.} This follows from a standard monodromy
argument, together with the analyticity of Hodge
subbundles.
\bigskip

Finally, we attend to: 

\bigskip
{\it Proof of claim.}\ It turns out that the singularities of $Y$ are quite mild. Note
that the singularities of $Y$ are introduced during the base change
$\Gamma\to \P^1$; $Y$ becomes singular when the map $\Gamma\to\P^1$ ramifies
over a point $t\in \P^1$ where $X_t$ is singular, i.e., it has an
ordinary double point. Therefore, the singularities of $Y$ have
the type of $x^2 + y^2 + z^2 + t^m = 0$. Let $p\in Y$ be such a singularity.
We may solve $p$ by a sequence of blowups:
\begin{equation}\label{E013}
\wt{Y} = Y_\mu \xrightarrow{\varphi_\mu}
Y_{\mu-1}\xrightarrow{\varphi_{\mu-1}} ... \xrightarrow{\varphi_1} Y_0 = Y
\end{equation}
where $\mu = \lfloor m/2 \rfloor$. The exceptional divisor \(E_k\subset
Y_k\) of $\varphi_k$ is a quadric in $\P^3$; it is a
cone over a conic curve if $2k < m$ and it is a smooth quadric if
$m = 2k$. Let $p_0 = p$ and $p_k\in E_k$ be the vertex of the cone
$E_k$ for $2k < m$. It is obvious that $Y_k$ is locally given by
$x^2 + y^2 + z^2 + t^{m - 2k} = 0$ at $p_k$ and \(\varphi_{k+1}:
Y_{k+1}\to Y_k\) is the blowup of $Y_k$ at $p_k$.

In order to pull back $\xi$ to $\wt{Y}$, we do it step by step, i.e.,
we first pull it back to $Y_1$, then $Y_2$ and so on. We will show
that there exists a sequence of cycles $\{\xi_k\in \CH^2(Y_k, 1)\}$
with all of them agreeing on the open set $Y\backslash \{p\}$.

By induction, it suffices to pull back the cycle \(\xi_{k-1}\in
\CH^2(Y_{k-1}, 1)\) to \(\xi_{k}\in
\CH^2(Y_{k}, 1)\).

Since $\varphi_k: Y_k\to Y_{k-1}$ is the blowup of $Y_{k-1}$ at
$p_{k-1}$,
\begin{equation}\label{E014}
Y_k\backslash E_k \isom Y_{k-1}\backslash \{p_{k-1}\}.
\end{equation}
So the question is again to extend a class in \(\CH^2(Y_k\backslash
E_k, 1)\) to $\CH^2(Y_k, 1)$. We look at the localization sequence
\begin{equation}\label{E015}
\CH^2(Y_k, 1) \xrightarrow{} \CH^2(Y_k\backslash E_k, 1)
\xrightarrow{}
\CH^1(E_k) \xrightarrow{\gamma} \CH^2(Y_k)
\end{equation}
If $E_k$ is a cone over a conic curve, then $\CH^1(E_k) = \BQ$ (see
\cite[Appendix A, Example 1.1.2, p. 428]{Ha}) and
$\gamma: \CH^1(E_k)\to \CH^2(Y_k)$ is obviously injective.

Suppose that $E_k$ is a smooth quadric. This happens in the last step
of blowups, i.e., when $k = \mu$ and $m = 2\mu$ is even.
Now
\begin{equation}\label{E016}
\CH^1(E_k) = \CH^1(\P^1\times \P^1) = \BQ\oplus \BQ.
\end{equation}
Let $L_1, L_2\subset E_k$ be the two rulings of $E_k$ which generate
$\CH^1(E_k)$. We claim that $L_1$ and $L_2$ are numerically
independent on $Y_k$, i.e., there exist divisors $D_1, D_2\subset Y_k$ such
that $D_i\cdot L_j = 0$ if $i = j$ and $D_i\cdot L_j \ne 0$ if $i\ne j$.
This certainly implies that $\gamma$ is injective.

Note that $Y_{k-1}$ has an ordinary double point \(x^2 + y^2 + z^2 +
t^2 = 0\) at $p_{k-1}$. It is well known that there exist two small
resolutions of $Y_{k-1}$. That is, we may blow down $Y_k$ along either
of the two rulings $L_1$ and $L_2$. Let $g: Y_k \to Y_k'$ be the
blowdown of $Y_k$ along $L_1$. Let $D$ be an ample divisor on
$Y_k'$. Then $g^* D \cdot L_2 \ne 0$ since $D$ is ample on $Y_k'$ and
$g^* D\cdot L_1 = 0$ since $g_* L_1 = 0$. We are done.

\end{document}